\documentclass[11pt]{amsart}
\usepackage{amsmath, amsthm, amssymb}
\usepackage{amssymb}
\usepackage{amsfonts}
\usepackage{amscd}
\usepackage{mathrsfs}
\usepackage{latexsym}
\usepackage{amsmath}
\usepackage{txfonts}
\usepackage{tikz}
\usepackage{tikz-cd}
\usepackage{graphicx}
\usepackage{subfig}
\usepackage{xcolor}
\usepackage{pdfsync}
\usepackage[all, cmtip]{xy}

\usepackage{mathtools}

\newtheorem{theorem}{Theorem}[section]
\newtheorem{lemma}[theorem]{Lemma}
\newtheorem{corollary}[theorem]{Corollary}
\newtheorem{proposition}[theorem]{Proposition}

\theoremstyle{definition}

\newtheorem{example}[theorem]{Example}

\theoremstyle{remark}

\newcommand{\spec}[1]{\textrm{spec}#1}
\newcommand{\sshf}[1]{\mathscr{O}_{#1}}
\newcommand{\shf}[1]{\mathscr{#1}}

\newcommand{\prj}[1]{\mathbb{P}^{#1}}

\newcommand{\iso}{\simeq}
\newcommand{\ses}[3]{0\rightarrow#1\rightarrow#2\rightarrow#3\rightarrow{0}}

\newcommand{\is}[1]{\mathscr{I}_{#1}}

\newcommand{\paren}[1]{\left(#1\right)}

\numberwithin{equation}{section}

\begin{document}
\allowdisplaybreaks
\title{Curves contracted by the Gauss map}


\author{Lei Song}

\address{School of Mathematics, Sun Yat-sen University, W. 135 Xingang Road, Guangzhou, Guangdong 510275, China}
\email{songlei3@mail.sysu.edu.cn}



\dedicatory{}


\begin{abstract}
Given a singular projective variety in some projective space, we characterize the smooth curves contracted by the Gauss map in terms of normal bundles. As a consequence, we show that if the variety is not linear, then a contracted line always has local obstruction  for the embedded deformation and each component of  the Hilbert scheme where the line lies is non-reduced everywhere.  
\end{abstract}

\maketitle

\section{Introduction}
Let $X\subseteq \prj{N}$ be a subvariety of dimension $n$. The Gauss map $\gamma: X\dashrightarrow \mathbb{G}(n, N)$ sends a smooth point $x\in X$ to the embedded tangent space $\widehat{T}_{x} X$ of $X$ at $x$.  

If $X$ is non-degenerate and smooth, then the map is a finite morphism by F. Zak \cite{Zak93}. For singular varieties, Griffiths-Harris \cite{GH79} established a structure theorem for Gauss map. Specifically it states that if  the dimension of the image $\dim (\text{Im}(\gamma))$ of $\gamma$ is less than $n$ (in this case, the Gauss map is called to be degenerate), then for a general point $[T]\in \text{Im}(\gamma)$, the closure of the fibre $\gamma^{-1}([T])$ is a projective space of dimension $n-\dim (\text{Im}(\gamma))$. Roughly speaking, a projective subvariety with the degenerate Gauss map can be built up by basic blocks, such as cones and developables. See  \cite{GH79}, \cite{AG01} for more precise statements, \cite{AG01} for a historical account, and \cite{Ran19} for a recent development. 

In general, it may happen that $\gamma$ is generically finite but has some positive-dimensional fibres, i.e.~it contacts some positive-dimensional subvarieties (see Example \ref{generically finite gp}). Here for a given subvariety $V\subseteq X$, we say $V$ is \textit{contracted} by the Gauss map $\gamma$ or say $\gamma$ is \textit{constant} on $V$  if $V$ is not contained in the singular locus $X_{\text{Sing}}$ and if the restricted map $\gamma|_V: V\dasharrow \mathbb{G}(n, N)$ is constant.

Projective subvarieties whose Gauss maps contract positive-dimensional subvarieties naturally appear in  several contexts, notably in A. Noma's classification for projective subvarieties (see \cite{Noma18} and the references therein). The usefulness of the classification is manifested in \cite{Noma14}, where a sharp bound for the Castelnuovo-Mumford regularity of $\sshf{X}$ was established for all non-degenerate, smooth varieties, and subsequently in \cite{MPS18} for an extension to isolated singularities. In the least specific class of Noma's classification, a projective subvariety has positive dimensional set of non-birational centers, which is a union of linear spaces (in characteristic 0), and the Gauss map is constant on each of the components.  It is worth further understanding of the geometry of such class of varieties. 

In this article, we study smooth curves that are contracted by the Gauss map. 
Our first main result is a characterization of contracted curves via the normal bundle.
\begin{theorem}\label{normal bundle formula}
Let $X\subseteq \prj{N}$ be a subvariety of dimension $n$ and $C\hookrightarrow X$ be a smooth projective curve not contained in $X_{\text{Sing}}$. Then the Gauss map is constant on $C$ if and only if the normal bundle $N_{C/X}$ of $C$ fits into the exact sequence
\begin{equation}\label{characterization sequence}
\ses{P^1(\sshf{C}(1))^*\otimes\sshf{C}(1)}{\oplus^{n+1} \sshf{C}(1)}{N_{C/X}}.
\end{equation}
In particular when $C=L$  is a line, the Gauss map is constant on $L$ if and only if 
\begin{equation}\label{formula for line}
  N_{L/X}\iso \oplus^{n-1} \sshf{L}(1).
\end{equation}
\end{theorem}

The sequence (\ref{characterization sequence}) is well-known when $X$ itself is a projective space, and the ``only if" part of the result may be viewed as a generalization.

Next let $L\subseteq X$ be a line contracted by the Gauss map. In view of the normal bundle formula (\ref{formula for line}),  one might expect at first glance that $L$ deforms in a family whose members cover $X$. Somewhat surprisingly, we show by dimension counting that as long as $X$ is not linear, $[L]$ \textit{always} has local obstruction for the embedded deformation in $X$. The details are presented in \S \ref{family}. As a consequence, we obtain
\begin{theorem}\label{non-reducedness}
Let $X\subseteq \prj{N}$ be a subvariety of dimension $n$ and $L\subseteq X$ be a line contracted by the Gauss map. Suppose $X\neq \prj{n}$. Then
any irreducible component $W$ of $\text{Hilb}^{o}$ containing $[L]$ is non-reduced everywhere, where $\text{Hilb}^{o}$ is the Hilbert scheme of lines on $X$ which are not contained in $X_{\text{Sing}}$. 
\end{theorem}

In \cite{Mumford62}, a component of the Hilbert scheme of space curves of degree 14 and genus 24 was proved to be nowhere reduced. Since then, more pathologies of Hilbert schemes have been found. Theorem \ref{non-reducedness} however presents pathologies of Hilbert schemes in a different flavor; here the parameterized objects are simply lines while the ambient space could be complicated. Example \ref{cone} shows that even very mild singularities, such as a cone over a smooth conic, can make the scheme non-reduced everywhere. We remark that  in case $X$ is normal, the result readily extends to rational curves of minimal degree by the nature of the proof and by \cite[Theorem 0.1]{CMS02}. 

In \S \ref{bound of degree}, we investigate upper bounds for the degree of a curve contracted by the Gauss map. Examples are given in \S \ref{example} to illustrate the results. 

Conventions and notations:  we work throughout over a base field $k$ which is algebraically closed and of characteristic 0. By a variety, we mean an integral, separated scheme, finite type over $k$. Given $k$-schemes $Z$ and $T$, we write $Z_T$ for $Z\times_{\spec{k}} T$. Given an embedding $i: X\hookrightarrow \prj{N}$, we put $V=H^0(\prj{N}, \sshf{\prj{N}}(1))$ and $\sshf{X}(1)=i^*\sshf{\prj{N}}(1)$. 

\subsection*{Acknowledgments}
I wish to take the opportunity to thank Xiaowen Hu, Joaqu{\'i}n Moraga, Jinhyung Park, and Ziv Ran for stimulating discussions. I am thankful to Yasuyuki Kachi and Nicholas Shepherd-Barron for answering my questions and clarifications.  I am also grateful to the anonymous referee for many valuable suggestions and for providing examples, which complement the paper. The work was partially supported by the Guangdong Basic and Applied Basic Research Foundation 2020A1515010876. 

\section{Preliminaries}
\subsection{Normal bundles of smooth curves}
We recall some general facts on the normal bundle of a smooth curve on a possibly singular variety.

Let $C$ be a smooth projective curve and $\shf{E}$ a coherent sheaf on $C$. Let $\tau$ be the torsion subsheaf of $\shf{E}$, and consider the exact sequence
\begin{equation*}
0\rightarrow \tau \rightarrow \shf{E} \rightarrow \shf{E} /{\tau} \rightarrow 0,
\end{equation*}
where $\shf{E} / \tau$ is locally free. Since $\text{Ext}^1(\shf{E} /{\tau}, \tau) \iso   H^1(C, \tau\otimes (\shf{E} /{\tau})^{\vee})=0 $ by the Grothendieck vanishing, the sequence splits
\begin{equation}\label{splitting of sheaves}
 \shf{E}\iso\tau\oplus\shf{E}/{\tau}.
\end{equation}

Let  $X$ be a projective scheme (possibly singular) and $C\hookrightarrow X$. The normal sheaf 
\begin{equation}\label{normal bundle}
  N_{C/X}\iso \shf{H}om_{\sshf{C}}(\is{C}/{\is{C}^2}, \sshf{C})\iso  \shf{H}om_{\sshf{C}}((\is{C}/{\is{C}}^2)/{\text{tor}}, \sshf{C})
\end{equation}
is locally free by (\ref{splitting of sheaves}). This justifies the terminology ``normal bundle".

 By the theory of $T^i$ functors (see \cite[Thm 3.5, Ex.~3.5, 3.6]{Hartshorne09}),
one has the exact sequence
\begin{equation*}
  0\rightarrow{T_C}\rightarrow{\shf{H}om(\Omega^1_X|_C, \sshf{C})}\rightarrow{N_{C/X}}\rightarrow \shf{T}^1(C, \sshf{C})\rightarrow \cdots,
\end{equation*}
where $T_C$ is the tangent bundle of $C$. Since $C$ is smooth, $\shf{T}^1(C, \sshf{C})=0$ cf.~ \cite[Ex.~4.3]{Hartshorne09}. Thus $N_{C/X}$ fits into the short exact sequence

\begin{equation}\label{short exact sequence for normal bundle}
  0\rightarrow{T_C}\rightarrow{\shf{H}om(\Omega^1_X|_C, \sshf{C})}\rightarrow{N_{C/X}}\rightarrow 0.
\end{equation}

\subsection{First order principle sheaf}
 Let $\shf{F}$ be a quasi-coherent sheaf on $X$. We briefly recall $P^1(\shf{F})$,  the first order principle sheaf of $\shf{F}$, with details referred  to the excellent exposition \cite[Chap. IV, Sec. A]{Kleiman77}. 
 
In general, one has the short exact sequence
\begin{equation}\label{extension for P^1}
\ses{\Omega^1_X\otimes\shf{F}}{P^1(\shf{F})}{\shf{F}}. 
\end{equation}
In complex analytic context, Atiyah \cite[Prop. 12]{Atiyah57} shows that if $\shf{F}$ is invertible, the extension class of (\ref{extension for P^1}) is indeed given by $-2\pi i c_1(\shf{F})\in H^1(X, \Omega^1_X)$.  By Lefschetz principle and Serre's GAGA, this result carries over to the algebraic setting immediately. 

When $X=\prj{1}$, an extension sheaf of $\sshf{\prj{1}}$ by $\Omega^1_{\prj{1}}$ is either trivial or isomorphic to $\sshf{\prj{1}}(-1)\oplus \sshf{\prj{1}}(-1)$. Thus for any $d\in \mathbb{Z}$, it holds 

\begin{equation}\label{P^1}
P^1(\sshf{\prj{1}}(d))\iso \left\{
\begin{array}{ll}
\sshf{\prj{1}}(d-1)\oplus \sshf{\prj{1}}(d-1) & d\neq 0 \\
\sshf{\prj{1}}\oplus \sshf{\prj{1}}(-2) & d=0. 
\end{array}
\right.
\end{equation}

\begin{lemma}\label{surjectivity}[\cite[(IV, 19)]{Kleiman77}]
For any closed embedding $X\hookrightarrow \mathbb{P}(V)$, the natural map
\begin{equation*}
  V\otimes \sshf{X}\twoheadrightarrow P^1(\sshf{X}(1))
\end{equation*}
is surjective.
\end{lemma}
This surjection induces an embedding $\mathbb{P}(P^1(\sshf{X}(1)))\hookrightarrow \mathbb{P}(V)\times X$ over $X$. In fact, over a smooth closed point $x\in X$, $\mathbb{P}(P^1(\sshf{X}(1))\otimes k(x))$ gives the embedded tangent space of $X$ at $x$.

\section{normal bundle and characterization}
In the section, we study the normal bundle $N_{C/X}$ of a smooth curve $C$ in a possibly singular subvariety $X\subseteq \prj{N}$. Put $\sshf{C}(1)=\sshf{X}(1)|_C$, so that the degree of a curve equals the degree of $\sshf{C}(1)$.  As usual, by a line, we mean a smooth rational curve of degree one. 

\begin{proposition}\label{normal bundle of line}
Let $C$ be a  smooth curve in $X$ of degree $d$. Then
$\shf{N}^*_{C/X}/\text{tor}\otimes\sshf{C}(d)$ is globally generated. 
Consequently when $C=L$ is a line,
  \begin{equation*}
   N_{L/X}\iso \oplus_i \sshf{L}(a_i),
  \end{equation*}
for some integers $a_i\le 1$.
\end{proposition}

\begin{proof}
Note to begin with that $\is{C/{\prj{N}}}(d)$ is globally generated, so the evaluation morphism
\begin{equation*}
  U\otimes_k\sshf{\prj{N}}(-d)\rightarrow \is{C/{\prj{N}}}
\end{equation*}
is surjective, where $U=H^0(\prj{N}, \is{C/{\prj{N}}}(d))$.
Pulling back to $X$, we have
\[U\otimes_k\sshf{X}(-d)\twoheadrightarrow  \is{C/{\prj{N}}}\otimes_{\sshf{\prj{N}}}\sshf{X}\iso \is{C/{\prj{N}}}/\paren{\is{X/{\prj{N}}}\cdot\is{C/{\prj{N}}}} \twoheadrightarrow \is{C/{\prj{N}}}/{\is{X/{\prj{N}}}}\iso \is{C/X},\]
which gives the surjection
\[
U\otimes_k \mathcal{O}_X(-d) \twoheadrightarrow \is{C/X}.
\]
Then pulling back to $C$, we get
\begin{equation*}
  U\otimes_k \sshf{C}(-d)\twoheadrightarrow \shf{N}^*_{C/X}.
\end{equation*}
Therefore
\begin{equation*}
 \shf{N}^*_{C/X}/\text{tor}\otimes\sshf{C}(d), 
\end{equation*}
as a quotient of $U\otimes_k \sshf{C}$ is globally generated.  When $C=L$ is a line, by taking dual, we get the assertion that $a_i\le 1$ for all $i$. 
\end{proof}

\begin{proposition}
Let  $C$ be a smooth curve not contained in $X_{\text{Sing}}$. Then there exists an exact sequence 
  \begin{equation*}
   \is{C/X}/{\is{C/X}^2}\otimes\sshf{C}(1)\rightarrow P^1(\sshf{X}(1))|_{C}\rightarrow P^1\left(\sshf{C}(1)\right)\rightarrow 0.
  \end{equation*}
 Consequently one has the short exact sequence
  \begin{equation}\label{exact sequence 1}
    \ses{\paren{\is{C/X}/{\is{C/X}^2}/{\text{tor}}}\otimes\sshf{C}(1)}{P^1(\sshf{X}(1))|_{C}/{\text{tor}}}{P^1\left(\sshf{C}(1)\right)}.
  \end{equation}
\end{proposition}
\begin{proof}
We have the commutative diagram with exact rows, see \cite[p.~343, (IV, 11)]{Kleiman77}
$$\xymatrix{
0\ar[r]& \Omega^1_X|_C\otimes\sshf{C}(1)\ar[r]\ar[d] & P^1(\sshf{X}(1))|_C\ar[r]\ar[d] & \sshf{C}(1)\ar[r]\ar[d]^{=} & 0\\
0\ar[r]& \Omega^1_C\otimes\sshf{C}(1)\ar[r]\ar[d] & P^1(\sshf{C}(1))\ar[r]\ar[d] & \sshf{C}(1)\ar[r]\ar[d] & 0\\
& 0 & 0 &0& .}$$
Moreover we have the natural exact sequence
\begin{equation*}
  \is{C/X}/{\is{C/X}^2}\rightarrow \Omega^1_X|_C\rightarrow\Omega^1_C\rightarrow 0.
\end{equation*}
Putting these together with the snake lemma, we obtain the exact sequence
\begin{equation*}
  \is{C/X}/{\is{C/X}^2}\otimes\sshf{C}(1)\rightarrow P^1(\sshf{X}(1))|_{C}\rightarrow P^1\left(\sshf{C}(1)\right)\rightarrow 0.
\end{equation*}
This yields the surjection of locally free sheaves
\begin{equation*}
  \paren{\is{C/X}/{\is{C/X}^2}/{\text{tor}}}\otimes\sshf{C}(1)\rightarrow
   \ker{\paren{P^1(\sshf{X}(1))|_{C}/{\text{tor}}\rightarrow P^1\left(\sshf{C}(1)\right)}}.
\end{equation*}
By the assumption that $C$ is not contained in the singular locus, the two bundles are of the same rank $n-1$, and hence are isomorphic.
\end{proof}

We need the lemmas in the sequel for Theorem \ref{normal bundle formula}. 

\begin{lemma}\label{global generation}
Let  $C$ be a smooth curve not contained in $X_{\text{Sing}}$. Then the vector bundle $P^1(\sshf{X}(1))|_{C}/{\text{tor}}$ is globally generated. 
\end{lemma}
\begin{proof}
In view of Lemma \ref{surjectivity}, we have the surjections
\[V\otimes\sshf{C}\twoheadrightarrow P^1(\sshf{X}(1))|_{C}\twoheadrightarrow P^1(\sshf{X}(1))|_{C}/{\text{tor}}. \]
Since $V\otimes\sshf{C}$ is globally generated, so is $P^1(\sshf{X}(1))|_{C}/{\text{tor}}$. 
\end{proof}

\begin{lemma}\label{criterion for trivial}
Let $C$ be a smooth projective curve and $E$ is a globally generated vector bundle of rank $r$. Then $E$ is trivial provided that either $h^0(C, E)=r$ or $\deg{E}=0$. 
\end{lemma}
\begin{proof}
By definition of globally generation, the evaluation map
\[\text{ev}: H^0(C, E)\otimes\sshf{C}\rightarrow E\]
is surjective. If $h^0(C, E)=r$, then the two bundles have the same rank. So ev is an isomorphism by Nakayama's lemma.  

Now assume $\deg{E}=0$. We shall proceed by induction on $r$. If $r=1$, a degree 0 line bundle with sections has to be trivial. For $r\ge 2$, take a nonzero section $s$ of $E$ and consider the exact sequence
\[0\rightarrow \sshf{C}\xrightarrow{s} E\rightarrow E'\rightarrow 0. \]
Suppose $E'$ has the torsion part $\tau$. Then by the snake lemma, we arrive at the commutative diagram with exact rows and columns
$$\xymatrix{
 &0\ar[d] &  & \tau\ar@{^(->}[d] & \\
0\ar[r]&\sshf{C}\ar[r]\ar[d] & E \ar[r]\ar@{=}[d] & E'\ar[r]\ar[d] & 0\\
0\ar[r]& L\ar[r]\ar@{->>}[d] & E \ar[r] & E'/{\tau}\ar[r] \ar[d]& 0.\\
& \tau &  &0 & }$$
Then $\deg{L}>0$, and hence $\deg(E'/{\tau})=-\deg{L}<0$. However this is impossible as $E'/{\tau}$ is globally generated. 

So  $E'$ is torsion free, and hence locally free. As a quotient of $E$, $E'$ is globally generated. By the induction hypothesis, $E'\iso \oplus^{r-1}\sshf{C}$.  It follows that 
\[r\le h^0(E)\le h^0(\sshf{C})+h^0(E')=r.\] 
This reduces to the case $h^0(C, E)=r$ and we finish the proof. 
\end{proof}

\begin{proposition}\label{criterion for constant Gauss}
Let $C\hookrightarrow X$ be a smooth curve not contained in $X_{\text{Sing}}$. Then the Gauss map is constant on $C$ if and only if
\[P^1(\sshf{X}(1))|_{C}/{\text{tor}}\iso \oplus^{n+1} \sshf{C}. \]
\end{proposition}

\begin{proof}
Put $E=P^1(\sshf{X}(1))|_{C}/{\text{tor}}$. Consider the commutative diagram
$$\xymatrix{
\mathbb{P}(E) \ar[rd]\ar@{^{(}->}[r]\ar@/^1.5pc/[rrrr]|f &  \mathbb{P}(P^1(\sshf{X}(1))|_{C})\ar[d]\ar@{^{(}->}[r] & \mathbb{P}(P^1(\sshf{X}(1))) \ar[d]\ar@{^{(}->}[r] & X\times \mathbb{P}(V)\ar[r]_{p_2}\ar[ld]^{p_1} & \mathbb{P}(V).  \\
& C\ar@{^{(}->}[r] & X  &  &
}$$
The first horizontal morphism on the top row is an injection because of the surjection $P^1(\sshf{X}(1))|_C \rightarrow E$,
the square in the middle is  Cartesian,
and the third horizontal morphism is an injection by Lemma~\ref{surjectivity}.

Denote by $f$ the composition of morphisms on the top row. The Gauss map on $C$ can be extended to a morphism from $C\rightarrow \mathbb{G}(n, N)$ in a unique way (see \cite[I.6.8]{Hartshorne77}), which corresponds to $f: \mathbb{P}(E)\rightarrow \mathbb{P}(V)$. In view of the diagram, $f$ is induced by a sub-linear system of $| \sshf{\mathbb{P}(E)}(1)|$.

Suppose the Gauss map is constant on $C$. Then the image of $f$ is the common embedded tangent space $\Lambda\iso \prj{n}$ of $X$ along $C$. Since $\sshf{\mathbb{P}(E)}(1)=f^*\sshf{\mathbb{P}(V)}(1)$ (or by Lemma \ref{global generation}), $E$ is globally generated. Moreover we have
\begin{equation*}
  h^0(C, E)=h^0(\mathbb{P}(E), \sshf{\mathbb{P}(E)}(1))=h^0(\Lambda, \sshf{\Lambda}(1))=n+1. 
\end{equation*}
So it follows from Lemma \ref{criterion for trivial} that $E$ must be trivial. 

Conversely if $E$ is trivial, then $h^0(\mathbb{P}(E), \sshf{\mathbb{P}(E)}(1))=h^0(C, E)=n+1$. So $f$ factors through an $n$-dimensional linear space $\Lambda\hookrightarrow \mathbb{P}(V)$, and hence the Gauss map is constant on $L$. This completes the proof. 
\end{proof}

\begin{proof}[Proof of Theorem \ref{normal bundle formula}]
Suppose the Gauss map $\gamma$ is constant on $C$. The claimed exact sequence is obtained by taking dual of (\ref{exact sequence 1}) and using Proposition \ref{criterion for constant Gauss}. 

Conversely,  by comparison  between sequence (\ref{characterization sequence}) and (\ref{exact sequence 1}), we deduce that the degree of the vector bundle $E=P^1(\sshf{X}(1))|_{C}/{\text{tor}}$ is 0.  On the other hand, $E$ is globally generated by Lemma \ref{global generation}. So it follows from Lemma \ref{criterion for trivial} that $E$ is trivial. Consequently $\gamma$ is constant on $C$ by Proposition \ref{criterion for constant Gauss}.

Consider the case $C=L$ be a line. In view of  Proposition \ref{normal bundle of line} and the fact $P^1\left(\sshf{L}(1)\right)\iso \oplus^2 \sshf{L}$, one sees that the sequence (\ref{exact sequence 1}) splits. 
Therefore $P^1(\sshf{X}(1))|_{L}/{\text{tor}}$ is trivial if and only if $N_{L/X}\iso \oplus^{n-1} \sshf{L}(1)$. This completes the proof.
\end{proof}

Observe that the normal bundle $N_{C/X}$ only depends on the embedding $C\hookrightarrow X$, but not on the embedding $X\hookrightarrow \prj{N}$, so the exact sequence (\ref{characterization sequence}) imposes strong restrictions on the possible situations where  $C$ is contracted by $\gamma$. As an example, we have

\begin{corollary}
Given $X\subseteq \prj{N}$ and $C\subseteq X$ be a smooth projective curve contracted by $\gamma$. Then for any ample line bundle $L$ on $X$ and a very ample linear system $U\subseteq H^0(X, \sshf{X}(1)\otimes L)$, the Gauss map $\gamma'$ associated to $X\subseteq \mathbb{P}(U)$ does not contract $C$.  In particular, any $d$-uple embedding does not contract $C$. 
\end{corollary}
\begin{proof}
Keep the notations as before. By Theorem \ref{normal bundle formula}, 
\[\deg N_{C/X}=(n+1)d_0+2g-2,\] 
where $d_0=\deg{\sshf{X}(1)|_C}$. But by the ampleness of $L$, $\deg(\sshf{X}(1)\otimes L)|_C>d_0$, so $N_{C/X}$ cannot fit into  (\ref{characterization sequence}). This shows $\gamma'$ does not contract $C$ by Theorem \ref{normal bundle formula}. 
\end{proof}

\section{Local obstructions and non-reducedness of Hilbert scheme}\label{family}
In this section, we study families of lines on a projective variety $X$ contracted by the Gauss map. 

Fix a very ample line bundle $\sshf{X}(1)$ and a polynomial $P\in\mathbb{Q}[t]$. By Grothendieck \cite{FGA}, the Hilbert scheme $\text{Hilb}^{P, \sshf{X}(1)}_{X/k}$ of $X$, whose closed points parameterize the closed subschemes of $X$ with the Hilbert polynomial $P$, exists and is projective. In particular for $P(t)=t+1$, we have
\begin{equation*}
  \text{Hilb}^{t+1, \sshf{X}(1)}_{X/k}\hookrightarrow \text{Hilb}^{t+1, \sshf{\prj{N}}(1)}_{\prj{N}/k}\iso \mathbb{G}(1, N)
\end{equation*}
parameterizing lines on $X$.

\begin{lemma}\label{large family of lines DNE}
Let $Y\subseteq \prj{N}$ be a subvariety of dimension $n$ and $y_0\in Y$ a closed point. Then there exists no $n$-dimensional family  of lines on $Y$ which pass $y_0$ and cover $Y$. 
\end{lemma}

\begin{proof}
Suppose not, then there exists an irreducible projective scheme $Z$ of $\dim n$ and  a flat family $\shf{U}$  of lines passing through $y_0$ over $Z$ such that the morphism $p$
$$\xymatrix{
\shf{U}\ar[d]\ar@{^{(}->}[r]\ar@/^1.5pc/[rr]|p & Z\times Y\ar[r]^{p_Y}\ar[ld]_{p_Z} & Y, \\
   Z                    &   &
}$$
is surjective. Thus for general $y\in Y$, it holds that $\dim p^{-1}(y)\ge 1$. This contradicts the fact that there exists at most one line passing through both $y_0$ and $y$. 
\end{proof} 

Let $Y$ be a $k$-scheme and $Z\hookrightarrow Y$ be a closed subscheme with defining ideal $\is{Z}$. Recall (cf. ~\cite{Kollar96}) the $\textit{obstruction space}$ $\text{Obs}(Z)\subseteq \text{Ext}^1_{\sshf{Y}}(\is{Z}, \sshf{Z})$ is the smallest $k$-subspace such that for each sequence $\ses{J}{B}{A}$, where $A, B$ are local Artin $k$-algebras with residue field $k$ and $J\frak{m}_B=0$, the obstruction class for extending a deformation of $Z$ over $A$ to one over $B$ lies in $\text{Obs}(Z)\otimes_k J$.

Now let $L\subset X$ be a line contracted by the Gauss map. The following result shows that in general, the corresponding closed point $[L]$ is a singular point of $\text{Hilb}^{1+t, \sshf{X}(1)}_{X/k}$.   

\begin{theorem}\label{obstruction}
Assume $X\neq\prj{n}$. Then $\text{Obs}(L)\neq (0)$. 
\end{theorem}
\begin{proof}
Let $W$ be an irreducible component containing $[L]$ with the universal family $\shf{U}$ above. Once again we have the following diagram
$$\xymatrix{
\shf{U}\ar[d]_{\pi}\ar@{^{(}->}[r]\ar@/^1.5pc/[rr]|p & W\times X\ar[r]^{p_X}\ar[ld]_{p_W} & X, \\
   W                     &   &
}$$
By \cite[(2.10.3)]{Kollar96} and Theorem \ref{normal bundle formula},
\[\dim W\ge h^0(L, N_{L/X})-\dim_k \text{Obs}(L)=2(n-1)-\dim_k \text{Obs}(L). 
\]
Since $W$ is proper, so is $p$. Thus $p(\shf{U})$ is closed in $X$. 

Suppose to the contrary that $\text{Obs}(L)=(0)$. Then $p$ has to be surjective, for otherwise there would exist $x_0\in X$, and irreducible $T\subseteq p^{-1}(x_0)$ with $\dim T\ge  n$. Then $\pi(T)$ parameterizes a family of lines on $p(\shf{U})$ passing through $x_0$ with $\dim T\ge n$ . This contradicts Lemma \ref{large family of lines DNE}.  

Now for each $x\in X$, $\dim p^{-1}(x)\ge n-1$. In view of Lemma \ref{large family of lines DNE},  we deduce that for each $x\in X$, 
$\dim p^{-1}(x)= n-1$.  Put $W_{x}=\pi(p^{-1}(x))$, which gives an $(n-1)$-dimensional family of lines passing through $x$, and consider the induced morphism $p_x: \shf{U}|_{W_x}\rightarrow X$. Because $\dim \shf{U}|_{W_x}=\dim X$ and there exists at most one line $l$ connecting $x$ and any other $x'$, we deduce that $p_x$ is surjective. This shows that for any $x\neq x'\in X$, there exists a line $l$ connecting $x$ and $x'$, implying that $X=\prj{n}$. This contradicts the assumption. 
\end{proof}

As a consequence we obtain:

\begin{corollary}
Assumptions as above. Then $L$ is locally obstructed. 
\end{corollary}
\begin{proof}
Suppose to the contrary that $L\subset X$ is locally unobstructed, then by \cite[I, Prop. 2.14]{Kollar96}, it holds that $\text{Obs}(L)\subseteq H^1(L, N_{L/X})$. By Theorem \ref{normal bundle formula}, $H^1(L, N_{L/X})=0$, which implies that $\text{Obs}(L)=(0)$, a contradiction to Theorem \ref{obstruction}. 
\end{proof}

\begin{lemma}\label{rigidity of constant Gauss map}
Given a flat family of smooth curves on $X$
$$\xymatrix{
\shf{C}\ar[rd]_{\pi}\ar@{^{(}->}[r] & X_T\ar[d]\ar@{-->}[r] & \mathbb{G}(n, N)_T\ar[ld], \\
                      & T  &
}$$
with the property that for each $t\in T$, $\shf{C}_t\nsubseteq (X_{\text{Sing}})_t$, where $X_{\text{Sing}}$ has the induced reduced scheme structure.  For each $t\in T$, one has the Gauss map $\gamma_{t}: \shf{C}_{t}\rightarrow \mathbb{G}(n, N)_t$. Assume the base scheme $T$ is irreducible and of finite type over $k$, and suppose that  for some closed point $t_0$, $\gamma_{t_0}$ is constant.  Then for general $t\in T$, $\gamma_t$ is constant. 
\end{lemma}
\begin{proof}
Note to begin with that given a curve $C\subseteq X$ and a change of the base field $\spec{K}\rightarrow \spec{k}$, the Gauss map $\gamma_K: C_K\dasharrow \mathbb{G}(n, N)_K$ is constant on $C_K$ if and only if $\gamma$ is constant on $C$.  Thus we can assume $T$ is reduced. Moreover, by taking a resolution of singularities $T'\rightarrow T$ and considering the family $\shf{C}_{T'}\rightarrow T'$,  we can assume that $T$ and hence $\shf{C}_{T}$ are smooth. By the valuative criterion for properness, one can extend the rational map $\shf{C}_{T}\dasharrow \mathbb{G}(n, N)_{T}$ over $T$ to a map defined over all codimension one points of $\shf{C}_T$. Therefore there exists a nonempty open subset $U\subseteq T$ over which we have the morphism $f:\shf{C}_{U}\rightarrow \mathbb{G}(n, N)_{U}$. 

Now let $\sshf{\mathbb{G}_U}(1)$ be the tautological line bundle on the Grassmannian $\mathbb{G}(n, N)_U$. It is evident that for each $t\in U$, $\gamma_t$ is constant if and only if $f^*\sshf{\mathbb{G}_U}(1)|_{\shf{C}_t}$ is trivial. Since $f^*\sshf{\mathbb{G}_U}(1)|_{\shf{C}_t}$ is globally generated, $f^*\sshf{\mathbb{G}_U}(1)|_{\shf{C}_t}$  is trivial if and only if its degree is zero. Then the assertion follows as the Hilbert polynomial is locally constant for fibres of a coherent sheaf flat over the base. 
\end{proof}

Since $Z=X_{\text{Sing}}\hookrightarrow X$, we have the closed embedding 
\[\text{Hilb}^{1+t, \sshf{X}(1)}_{Z/k}\hookrightarrow \text{Hilb}^{1+t, \sshf{X}(1)}_{X/k}.\] Put $\text{Hilb}^{o}=\text{Hilb}^{1+t, \sshf{X}(1)}_{X/k}\backslash \text{Hilb}^{1+t, \sshf{X}(1)}_{Z/k}$ with the induced scheme structure, and let $W$ be an irreducible component of $\text{Hilb}^{o}$.

\begin{proof}[Proof of Theorem \ref{non-reducedness}]
By Lemma \ref{rigidity of constant Gauss map}, for a general point $x: \spec k\rightarrow W$, the Gauss map  on the corresponding line $L_x$ is constant, and hence $x$ is singular by Theorem \ref{obstruction}.  Suppose to the contrary that $W$ is reduced at some point $x'$, then there exists an open neighborhood $U\subseteq W$ of $x'$ such that $U$ is reduced with the induced scheme structure.  Then by the generic smoothness, there exists a nonempty open $V\subseteq U$ that is smooth over $k$. This gives a contradiction.
\end{proof}

\section{Bound of degree}\label{bound of degree}
Gauss maps which are generically finite but not quasi-finite seem much less understood than the degenerate ones.  
Mezzetti-Tommasi give an example (cf.~\cite[Remark 2]{Ran19}) where a singular quartic surface has a generically finite Gauss map which contracts a smooth plane quadric. 

In surface case, we give an upper bound of the degree of the contracted curves. Note in higher dimensions, such a bound cannot exist because the closure of a fibre of the Gauss map can be $\prj{m}$ with $m\ge 2$. 

\begin{proposition}\label{rational curve of high degree is not contracted}
Let $X$ be a non-degenerate singular surface in $\prj{N}$. Suppose $\gamma$ is constant on a smooth projective curve $C$. Then $C$ is a plane curve and 
\[\deg{X}\ge 2\deg{C}+(N-3).\] 
\end{proposition}

\begin{proof}
It is clear that $C$ is a plane curve contained in the embedded tangent space $\widehat{T}=\widehat{T}_xX$ for general $x\in C$. 
When $N=3$, since $X$ and $\widehat{T}$ intersect at $C$ doubly, by B\'ezout's theorem, we have $\deg X\ge 2\deg C$.  For $N>3$, note the secant variety $\text{Sec}(C)\subseteq \widehat{T}$. Let $\overline{C}(X)$ be the closure of the non-birational centers of $X$, which is a proper linear space in $X$. Then for any smooth point $x_0\in X\backslash (\text{Sec}(C)\cup \overline{C}(X))$,  the inner projection $\pi: X\backslash\{x_0\}\rightarrow \prj{N-1}$ gives rise to a birational map $X\dasharrow X':=\overline{\pi(X\backslash\{x_0\})}$ such that the restriction $\pi|_C: C\rightarrow C':=\pi(C)$ is an isomorphism. So $X'\subseteq \prj{N-1}$ is a surface with a smooth rational curve $C'$ on it, with $\deg X'=\deg X-1$ and $\deg C'=\deg C$.  Fix a closed point $y_0\in C\backslash X_{\text{Sing}}$. By taking the center $x_0$ general, we can guarantee that $\pi$ is isomorphic at $y_0$, and hence is isomorphic at a neighborhood $U$ of $y_0$.  Therefore for $y\in U$, $\widehat{T}_{\pi(y)}(X')=\pi(\widehat{T}_y(X))=\pi(\widehat{T})$, which implies that $\gamma_{X'}$ is constant on $C'$ provided $x_0$ is taken general. So we can proceed by induction on $N$. 
\end{proof}

In a different direction, we establish an upper bound for the degree of contracted curves by the Gauss map with respect to an anticanonical divisor. 

\begin{theorem}\label{anticanonical degree}
Let $X\subseteq \prj{N}$ be a normal $\mathbb{Q}$-Gorenstein variety of dimension $n$. Let $C\subseteq X$ be a smooth projective curve of degree $d$.  Suppose the Gauss map is constant on $C$. Then $-K_X\cdot C\le (n+1)d$. 
\end{theorem}
\begin{proof}
Suppose $rK_X$ is Cartier for some integer $r>0$. Put $\Omega^n_X=\wedge^n \Omega_X$. We have the natural morphism $\eta: \Omega^n_X\rightarrow \omega_X$ (cf. ~\cite[Appendix]{EM09}), where $\omega_X$ is the canonical sheaf.  The morphism $\eta$ induces $(\Omega^n_X)^{[r]}\rightarrow \omega^{[r]}_X\iso \sshf{X}(rK_X)$, where $\star^{[r]}$ denotes the double dual of $\star^{\otimes r}$. Consider the composite of the natural morphisms
\begin{equation*}
\xi: ((\Omega^n_X)|_C)^{\otimes r}\iso (\Omega^n_X)^{\otimes r}|_C \rightarrow (\Omega^n_X)^{[r]}|_C\rightarrow \sshf{X}(rK_X)|_C.
\end{equation*}
Since $\xi$ is generically injective, the induced morphism 
\[((\Omega^n_X)|_C)^{\otimes r}/{\text{tor}}\rightarrow \sshf{X}(rK_X)|_C\]
is injective. Consequently
\[rK_X\cdot C=\deg \sshf{X}(rK_X)|_C\ge \deg ((\Omega^n_X)|_C)^{\otimes r}/{\text{tor}}. \]
To finish the proof,  note by (\ref{short exact sequence for normal bundle}) and Theorem \ref{normal bundle formula}, $\shf{H}om(\Omega^1_X|_C, \sshf{C})$ is a vector bundle of degree $(n+1)d$, so
$(\Omega^n_X|_C)/{\text{tor}}$ is a line bundle of degree $-(n+1)d$. 
\end{proof}

\section{Examples}\label{example}
In this section, keeping notation as before, we give several examples to illustrate results in the previous sections. 
\begin{example}\label{cone}
Let $Y$ be a non-degenerate smooth projective variety dimension $n-1$ in a hyperplane $H$ of $\prj{N}$. Let $X$ be the cone over $Y$ with the vertex $v\notin H$. The singular locus of $X$ is $v$. For any $x\neq v$, there exists a unique point $y\in Y$ such that $x\in L= \langle{v, x}\rangle$.  By Terracini's lemma, we deduce that
\[\langle {v, \widehat{T}_{y} Y}\rangle = \widehat{T}_{x} X. \]
In particular, the ruling $L$ is contracted by the Gauss map.  These rulings are naturally parameterized by $Y$ as a topological space. Precisely, we have
\[Y=(\text{Hilb}^0X)_{\text{red}}. \]

To be concrete and for simplicity, take $Y$ to be the smooth conic on the hyperplane $x_3=0$ in $\prj{3}$, defined by $x_0x_1-{x_2}^2=0$, and take $v=[0:0:0:1]$. Then the cone $X$ over $Y$ with vertex $v$ has only canonical singularities.  

For $I=\{2, 3\}\subset\{1, 2, 3, 4\}$, the points of the affine open $U_I$ of $G(2, 4)\iso \mathbb{G}(1, 3)$ are parameterized by the matrices of the form
 \[
\begin{bmatrix}
1 & a& b & 0\\
0 & c & d & 1\\
\end{bmatrix}
\]
Then 
\[\text{Hilb}^0X\cap U_{\{2, 3\}}\iso \spec k[a, b, c, d]/(a-b^2, d^2, c-2bd)\iso \spec k\left[b, d\right]/(d^2).\]
See Appendix for an elaboration. Similarly, in $U_{\{1, 3\}}$ whose points are parametrized by the matrices of form 
 \[
\begin{bmatrix}
a & 1& b & 0\\
c & 0 & d & 1\\
\end{bmatrix},
\]
it also holds that
\[\text{Hilb}^0X\cap U_{\{1, 3\}}\iso \spec k\left[b, d\right]/(d^2). \]
Since $\text{Hilb}^0X \subseteq U_{\{2, 3\}}\cup U_{\{1, 3\}}$, we see that $\text{Hilb}^0X$ is non-reduced everywhere. 

Any line on the cone $X$ over a conic has to pass the vertex, so  in this example, we actually have:
\[Y=(\text{Hilb}^{t+1, \sshf{X}(1)}_{X/k})_{\text{red}.}\]
\end{example}

Let $X$ be a singular hypersurface defined by $F$ in $\prj{N}$. The singular locus $X_{\text{Sing}}$ has a natural scheme structure given by the Jacobian ideal  $(\frac{\partial{F}}{\partial x_0}, \cdots, \frac{\partial{F}}{\partial x_N})$. With the scheme structure, $X_{sing}$ is called the singular scheme of $X$, and there are strong constraints on $X_{\text{Sing}}$, as shown in \cite{Aluffi95}.  In the following example, among all lines, the ones contracted by the Gauss map of $X$ have maximal intersection with $X_{\text{Sing}}$. 

\begin{example}
Consider the Veronese surface $\nu: \prj{2}\hookrightarrow \prj{5}$, $[x_0: x_1: x_2]\mapsto [x^2_0: x_0x_1: x_0x_2: x^2_1: x_1x_2: x^2_2]$. The ideal $I_{\nu(\prj{2})}$ is generated by $2\times 2$ minors of 
\[
\begin{bmatrix}
T_0 & T_1 & T_2\\
T_1 & T_3 & T_4\\
T_2 & T_4 & T_5
\end{bmatrix}
\]
The determinant
\[F=T_0T_3T_5-T_0T^2_4-T^2_1T_5-T^2_2T_3+2T_1T_2T_4\]
defines a cubic hypersurface $X$. The Jacobian ideal of $X$ is precisely $I_{\nu(\prj{2})}$; $\nu(\prj{2})$ is the singular scheme of $X$ in the terminology of \cite{Aluffi95}. 

Given any line $l\subset \prj{2}$, the smooth conic $C=\nu(l)$ is contained in a unique plane $\Pi\iso \prj{2}$.  For any line $L\subset \Pi$, $L\subset X$. Thus $\Pi\subset X$. In fact, $\Pi$ is contracted by the Gauss map $\gamma_X: X\backslash \nu(\prj{2})\rightarrow \mathbb{G}(4, 5)$. 

To see this, we can assume $l: x_2=0$. Then $\Pi$ is given by $T_2=T_4=T_5=0$, and $C$ is cut out by the extra equation $T^2_1-T_0T_3=0$.  For any  point $P=[p_0:\cdots : p_5]\in 
\Pi\backslash C$, 
\[\nabla F(P)=\left[0, 0, 0, 0, 0, p^2_1-p_0p_3\right],\]
indicating that the projective tangent space at $P$ is given by $T_5=0$, independent of $P$. 

Furthermore, there exists no 3 dimensional linear space in $X$, so $\gamma_X$ is equi-dimensional. 

The Hilbert scheme of lines $\text{Hilb}^{t+1, \sshf{X}(1)}_{X/k}$ has two irreducible components $W_1$ and $W_2$, both of dimension 4. The general point of $W_1$ is a line which lies on an embedded tangent space to $\nu(\prj{2})$, and which does not intersect with $\nu(\prj{2})$, while the general point of $W_2$ is a line which lies on a $\Pi$. By standard computation, for a general $[L]\in W_1$, the normal bundle 
\[N_{L/X}\iso \sshf{}(1)\oplus \sshf{}(1) \oplus \sshf{}(-1), \]
implying that $W_1$ is smooth at $[L]$.  On the other hand, by (\ref{formula for line}) of Theorem \ref{normal bundle formula}, for any line $[L]\in W_2$, 
\[N_{L/X}\iso \sshf{}(1)\oplus \sshf{}(1)\oplus \sshf{}(1),\]
which implies $W_2$ is non-reduced everywhere. 

\end{example}

The following example is due to the referee. It gives a recipe to construct a generically finite Gauss map that contracts subvarieties of positive dimension.
\begin{example}\label{generically finite gp}
Let $\tilde{X}\subseteq \prj{N}$ be a dual defective, smooth subvariety of dimension $n$. Note that $\tilde{X}$ cannot be a hypersurface. Denote the dual defect by $\delta>0$. By the generic smoothness, there exist a smooth subvariety $\tilde{Y}$ of dimension $\delta$ and a hyperplane $H$ such that $H$ is tangent to $\tilde{X}$ precisely along $Y$. Now take a general linear subspace $\Lambda\iso \prj{N-n-2}$ of $H$, and consider the projection  $\pi: \prj{N}\dashrightarrow \prj{n+1}$ from the center $\Lambda$. Let $X=\pi(\tilde{X})$ and $Y=\pi(\tilde{Y})$. 

For general $x\in Y$,  $\widehat{T}_{x} X=H\cap \prj{n+1}$, while for $x\notin Y$ and $x$ smooth, $\widehat{T}_{x} X$ is different from $H\cap \prj{n+1}$. Therefore $Y$ is a (topological) fibre of the Gauss map $\gamma_X$ of $X$. Since $Y$ is not linear, it follows that $\gamma_X$ is generically finite. 

Note if $(2\delta+1)+(N-n-2)<N-1$, i.e. $\delta<\frac{n}{2}$, one can arrange $\Lambda$ so that $Y\iso \tilde{Y}$, so make $Y$ smooth. 
\end{example}

\section*{Appendix: Local equations for Hilbert schemes of lines}
In the appendix, we explain given a projective variety $X\hookrightarrow \prj{N}$, how to write down the local equations for $\text{Hilb}^{t+1, \sshf{X}(1)}_{X/k}$ in $\mathbb{G}(1, N)$.  Recall
\[ \mathbb{G}(1, N)\iso G(2, N+1)=\bigcup_I U_I,\]
where $I$ ranges over subsets of $\{1, 2, \cdots, N+1\}$ with $|I|=2$, and $U_I\iso \mathbb{A}^{2(N-1)}$, whose closed points are represented by $2\times (N+1)$ matrices with $I$-th minors being the $2\times 2$ identity matrix.

For each $I$, we have the Cartesian diagram
$$\xymatrix{
\text{Hilb}^{t+1, \sshf{X}(1)}_{X/k, I} :=\text{Hilb}^{t+1, \sshf{X}(1)}_{X/k}\cap U_I \ar[d] \ar@{^{(}->}[r] & U_I\ar[d]\\
\text{Hilb}^{t+1, \sshf{X}(1)}_{X/k} \ar@{^{(}->}[r] &  G(2, N+1).} $$
To describe the defining equations for $\text{Hilb}^{t+1, \sshf{X}(1)}_{X/k, I}$ in $U_I$, we need the simple lemmas.

\begin{lemma}\label{local triviality}
Let $T$ be a locally Noetherian $k$-scheme. Let $\shf{L}\hookrightarrow \prj{N}_T$ be a flat family of lines over $T$. Denote by $\pi$ the morphism $\shf{L}\rightarrow T$. Then $\shf{L}\iso \mathbb{P}(\pi_*\sshf{\shf{L}}(1))$, i.e. $\shf{L}$ is a projective bundle over $T$. 
\end{lemma}
\begin{proof}
By flat base change, 
\[\pi_*\sshf{\shf{L}}(1)\otimes_{\sshf{T}}k(t)\xrightarrow{\iso} H^0\left(\shf{L}_t, \sshf{\shf{L}_t}(1)\right), \]
implying that $\pi_*\sshf{\shf{L}}(1)$ is a rank 2 vector bundle on $T$. 
On $\shf{L}$, we have the natural morphism
\[\pi^*\pi_*\sshf{\shf{L}}(1)\rightarrow \sshf{\shf{L}}(1),\]
which is surjective. As a result, we have
\[\shf{L}  \xrightarrow{\varphi}  \mathbb{P}\left(\pi_*\sshf{\shf{L}}(1)\right)\]
 over $T$ with the property that
 \[\varphi^* \sshf{\mathbb{P}(\pi_*\sshf{\shf{L}}(1))}(1)\iso \sshf{\shf{L}}(1). \]
Thus $\varphi$ induces an isomorphism. 
\end{proof}

\begin{lemma}\label{local check}
Let $f: Y\rightarrow X$ be a morphism of schemes and $i: Z\hookrightarrow X$ a closed subscheme. Let $\{V_j\}$ be an open covering of $Y$ and let $f_j$ be the restriction of $f$ to $V_j$. Then there exists a morphism $g: Y\rightarrow Z$ such that $i\circ g=f$ if and only if for each $j$, there exists a morphism  $g_j: V_j\rightarrow Z$ such that $i\circ g_j=f_j$. 
\end{lemma}
\begin{proof}
The existence of $g\iff \is{Z/X}\cdot \sshf{Y}=0 \iff \is{Z/X}\cdot \sshf{V_j}=0$ for all $j\iff$ the existence of $g_j$ for all $j$. 
\end{proof}

Let $T$ be a locally Noetherian $k$-scheme,
\begin{eqnarray*}
&& Hom_k(T, U_I)\\
&\iso &\{\text{$T$-points of  $G(2, N+1)$ that factor through } U_I\}\\
&\iso & \{ 
\shf{L}  \hookrightarrow \prj{N}_T \:|\: \shf{L} \text{ is flat over } T; \forall t\in T, \shf{L}_t\hookrightarrow \prj{N}_t\text{ is a line of type I over } k(t)\}.
\end{eqnarray*}

Given $a\in Hom_k(T, U_I)$ and let $\shf{L}$ be the associated family of lines over $T$. Then $a\in Hom_k(T, \text{Hilb}^{t+1, \sshf{X}(1)}_{X/k, I})$ if and only if 
\[\xymatrix{
\shf{L} \ar@{-->}[rd] \ar@{^{(}->}[rr] & & \prj{N}_T\\
 &  X_T\ar@{^{(}->}[ru] & }\]
and by Lemma \ref{local check}, if and only if
\[\xymatrix{
\shf{L}_{T_j} \ar@{-->}[rd] \ar@{^{(}->}[rr] & & \prj{N}_{T_j}\\
 &  X_{T_j}\ar@{^{(}->}[ru] & }\]
 for any open covering $\{T_j\}$ of $T$.

In view of Lemma \ref{local triviality}, there exists an open covering $\{T_j\}$ of $T$ such that $\shf{L}_{T_j}$ are trivial for all $j$.  We assume that $\shf{L}\iso \text{Proj}\sshf{T}[s, t]$, then the preceding factorization amounts to one in the diagram of $\sshf{T}$-algebras:
\[\xymatrix{
 & \sshf{T}[s, t] \\
\sshf{T}[x_0, x_1, \cdots, x_N]/{I_X} \ar@{-->}[ru]^{\bar{\psi}}&  \sshf{T}[x_0, x_1, \cdots, x_N], \ar[u]_{\psi}\ar[l] }
\]
where $I_X$ is the saturated homogeneous ideal defining $X$, and $\psi(x_i)\in H^0\left(\mathbb{P}^1_T, \sshf{\mathbb{P}^1_T}(1)\right)$. 

The existence of such $\bar{\psi}$ in turn is equivalent to 
\begin{equation}\label{factorization criterion}
\psi(I_X)=(0).
\end{equation}
 
Keeping in mind that $a\in Hom_k(T, U_I)$, $\psi(x_i)$ can be described explicitly as follows. For ease of notation,  assume that $I=\{1, 2\}$. Lines in the family $\shf{L}$ are represented by the matrices
\[\begin{bmatrix}
1 & 0 & a_{1, 1} & \cdots & a_{1, N-1}\\
0 & 1 & a_{2, 1} & \cdots & a_{2, N-1}
\end{bmatrix},
\]
where $a_{ij}\in \Gamma(T, \sshf{T})$ are determined by $a$. As a result,
\begin{equation}\label{describing lines}
\psi(x_0)=s, \psi(x_1)=t, \psi(x_3)=a_{1, 1}s+a_{2, 1}t, \cdots, \psi(x_{N})=a_{1, N-1}s+a_{2, N-1}t. 
\end{equation}
Plugging (\ref{describing lines}) into (\ref{factorization criterion}), we obtain functions: 
\[h_1, \cdots, h_m\in A(U_I)= k\left[x_{1,1}, \cdots, x_{1, N-1}; x_{2,1}, \cdots, x_{2, N-1}\right],\]
in the affine coordinate ring of $U_I$, such that
\[h_l(a_{i, j})=0 \quad\qquad \text{ in } \:\Gamma(T, \sshf{T})\]
for all $1\le l\le m$. Therefore (\ref{factorization criterion}) is equivalent to 
\[\xymatrix{
 &\Gamma(T, \sshf{T}) \\
A(U_I)/(h_1, \cdots, h_m) \ar@{-->}[ru] &  A(U_I). \ar[u]_{a^*}\ar[l] }
\]
In summary, we have obtained
\begin{proposition}
In $U_I$, $\text{Hilb}^{t+1, \sshf{X}(1)}_{X/k, I}$ is defined scheme-theoretically by $h_1, \cdots, h_m$.
\end{proposition}
 
\bibliography{Gaussmap.bib}{}
\bibliographystyle{alpha}

\end{document}